\documentclass[12pt,english]{article}
\usepackage{babel}
\usepackage[cp1251]{inputenc}
\usepackage{amsmath,amssymb}
\begin{document}

\newpage

\begin{center}
{\bf \Large
Localization of the matrix spectrum \\ and Lyapunov type equations\footnote{The study was carried out
within the framework of the state contract of the Sobolev Institute of Mathematics 
(project no.~FWNF-2022-0008).}}
\end{center}

\begin{center}
{\bf G.V.~Demidenko$^{1,2}$, Z.~Wang$^{2}$}
\end{center}

\begin{center}
${}^1$Sobolev Institute of Mathematics, \\
${}^2$Novosibirsk State University, \\
Novosibirsk, Russia
\end{center}

\begin{center}
demidenk@math.nsc.ru
\end{center}

{\bf Abstract.}
The problems on the location of the matrix spectrum inside or outside domains
bounded by ellipses or parabolas are studied.
Special Lyapunov-type equations are connected with these problems.
Theorems about the unique solvability of such equations are proved.
Conditions for perturbations of matrix entries are obtained,
which guarantee that the spectra of the perturbed matrices belong to the specified domains 
as well.

\vskip10pt


\vskip10pt


\section{Introduction}

In the paper we consider matrix equations of the form
$$
\sum\limits^N_{j,k=0} a_{jk} (A^*)^jHA^k = C,
\eqno(1)
$$
where
$A$, $C$, 
are given 
$(n \times n)$-matrices,
$a_{jk}$
are numerical coefficients.
We study the relationship between the solvability of such equations
and the belonging of the spectrum of 
$A$
to sets lying inside or outside domains bounded by an ellipse or a parabola.
Conditions for perturbations of matrix entries are obtained, which guarantee
that the spectrum of the perturbed matrix
$A + B$
belongs to the specified domains as well.

Note that many problems of control theory, problems on description of
$\varepsilon$-spectrum
or pseudospectrum lead to the necessity to study the belonging
of the matrix spectrum to various domains bounded
by some contours in the complex plane, as well as to solve
the problem on dichotomy of the matrix spectrum with respect to some curves.
Therefore, the establishment of various criteria and the development of algorithms
for determining the location of the matrix spectrum in the complex plane are of great interest.

The most well-known criteria on the location of the matrix spectrum are the Lyapunov criteria
(see, for example,~\cite{DK},~\cite{G}).
These criteria provide necessary and sufficient conditions 
under which the matrix spectrum belongs to the left half-plane
$$
{\cal C}_- = \{\lambda \in {\mathbb C}: \ {\rm Re}\, \lambda < 0\}  
$$ 
or the unit disk
$$
{\cal C}_i = \{\lambda \in {\mathbb C}: \ |\lambda| < 1\}. 
$$ 
We give these criteria below.

\vskip5pt

{\bf Theorem 1.}
{\it 
The spectrum of the matrix
$A$
belongs to the left half-plane
${\cal C}_-$ 
if and only if the equation
$$
HA + A^*H = -C 
\eqno(2) 
$$
with 
$C = C^* > 0$
has a Hermitian positive definite solution
$H = H^* > 0$.
}

\vskip5pt

In the literature, (2) is called the {\it Lyapunov equation}.

Recall that (2) arises when studying the stability of solutions
to systems of ordinary differential equations of the form
$$
\frac{dy}{dt} = Ay + f(y).
\eqno(3)
$$

Another interesting example is
the {\it discrete Lyapunov equation}
$$
H - A^*HA = C,
\eqno(4)
$$
where
$A$, $C$
are given 
$(n \times n)$-matrices,
$H$
is an unknown matrix.

\vskip5pt

{\bf Theorem 2.}
{\it
The spectrum of the matrix
$A$
belongs to the unit disk
${\cal C}_i$ 
if and only if
$(4)$ with 
$C = C^* > 0$
has a Hermitian positive definite solution
$H=H^* > 0$.
} 

\vskip5pt

As known, (4) arises when studying the stability
of solutions to difference equations of the form
$$
y_{j+1} = Ay_j + \varphi_j(y),
\qquad
j = 0,1,\dots.
\eqno(5)
$$

Obviously, (2) and~(4) belong to the class of equations of the form~(1).
Many examples of equations of this type are given in the monograph~\cite{M}
when solving problems on location of the matrix spectrum in different domains,
boundaries of which are defined by algebraic equations.
Note that the study of such problems is based on M.G.~Krein's theorem
on the solvability of a class of {\it generalized Lyapunov equations}
of the following form
$$
\sum\limits^N_{j,k=0} a_{jk} B^jHA^k = Y,
\eqno(6)
$$
where
$A$, $B$, $Y$
are given matrices of sizes
$n \times n$, $m \times m$, $m \times n$,
respectively,
$N$
defines the order of (6),
$a_{jk}$
are given constant coefficients, and
$H$
is an unknown matrix of size
$m \times n$. 
We formulate this theorem below~\cite{DK}.

Introduce the following notation:
$$
P(\lambda,\mu) = \sum\limits^N_{j,k=0} a_{jk}\lambda^j \mu^k,
$$
where
$a_{jk}$
are the coefficients in (6);
$\sigma(A) = \{\mu_1,\dots,\mu_n\}$
is the spectrum of the matrix
$A$; 
$\sigma(B) = \{\lambda_1,\dots,\lambda_m\}$
is the spectrum of the matrix
$B$; 
$\gamma_A$
is a contour in the complex plane surrounding 
$\sigma(A)$; 
$\gamma_B$
is a contour surrounding 
$\sigma(B)$.

\vskip5pt

{\bf Theorem 3} (M.G.~Krein~\cite{DK}).
{\it
Let 
$$
P(\lambda_s,\mu_r) \neq 0,
$$
where
$$
\lambda_s \in \sigma(B), \quad s = 1, \dots,m,
\qquad
\mu_r \in \sigma(A), \quad r =1,\dots,n.
$$
Then there is a unique solution to $(6)$ for any matrix
$Y$
and it has the form
$$
H = \frac{1}{(2\pi i)^2} \int\limits_{\gamma_A} \int\limits_{\gamma_B}
\frac{1}{P(\lambda, \mu)} (\lambda I - B)^{-1} Y (\mu I - A)^{-1}
d\lambda\, d\mu.
$$
}

\vskip5pt

Various conditions for the belonging of the matrix spectrum to bounded domains,
which are formulated in terms of the solvability of matrix equations of the form~(6),
can be found in~\cite{DK}--\cite{D09}.

In this paper, we consider two matrix equations of the form~(1)
that arise when studying the localization of the matrix spectrum
in the complex plane with respect to an ellipse and a parabola.
Solvability criteria are known for these equations (see~\cite{M},~\cite{D09}).
Our aim is to consider these equations with perturbations of matrix entries,
to obtain conditions for the perturbations under which the localization
of the matrix spectrum does not change, and to study the solvability of these equations
with perturbed matrix entries.

\section{Localization of the matrix spectrum with respect to an ellipse}

The following problem was considered in~\cite{M},~\cite{D09},~\cite{D01}.

\vskip5pt

{\bf Problem 1.} Belonging of the spectrum of the matrix
$A$
to the domain bounded by an ellipse
$$
{\cal E}_i = 
\left\{\lambda : \ \frac{({\rm Re}\,\lambda)^2}{a^2}
+ \frac{({\rm Im}\,\lambda)^2}{b^2} < 1 \right\},
\quad
a > b,
$$
is equivalent to the solvability of the matrix equation
$$
H - \left(\frac{1}{2a^2}+\frac{1}{2b^2}\right) A^*HA 
- \left(\frac{1}{4a^2}-\frac{1}{4b^2}\right)(HA^2+(A^*)^2H) = C
\eqno(7)
$$
with 
$C = C^* > 0$,
in the class of Hermitian positive definite matrices
$H=H^*>0$. 

\vskip5pt

In particular, the following theorem was proved in~\cite{D01}.

\vskip5pt

{\bf Theorem 4.}
{\it
Let
$C = C^* > 0$. 
The spectrum of the matrix
$A$
belongs to 
${\cal E}_i$ 
if and only if there is a Hermitian positive definite solution
$H = H^* > 0$ 
to~$(7)$.
}

\vskip5pt

The study of some problems related with~(7)
was carried out in~\cite{D01},~\cite{BDM}.

Note that (7) has the form (6) with 
$$
N=2, \quad n=m, \quad B=A^*, \quad Y=I,
$$
$$
a_{00}=1, \quad a_{01}=a_{10}=a_{12}=a_{21}=a_{22}=0, 
$$
$$ 
a_{11}=-\left(\frac{1}{2a^2}+\frac{1}{2b^2}\right), \quad 
a_{02}=a_{20}=-\left(\frac{1}{4a^2}-\frac{1}{4b^2}\right).
$$

Consider the problem on the belonging of the matrix spectrum to a set
lying outside the closure of the domain bounded by an ellipse,
i.e. the problem on the belonging to the set
$$
{\cal E}_e = 
\left\{\lambda : \ \frac{({\rm Re}\,\lambda)^2}{a^2}
+ \frac{({\rm Im}\,\lambda)^2}{b^2} > 1 \right\},
\quad
a > b.
$$

By analogy with the problem on the belonging of the matrix spectrum to 
${\cal E}_i$,
consider the matrix equation of the form~(7) with a modified right-hand side
$$
H - \left(\frac{1}{2a^2}+\frac{1}{2b^2}\right) A^*HA 
- \left(\frac{1}{4a^2}-\frac{1}{4b^2}\right)(HA^2+(A^*)^2H) = - C.  
\eqno(8)
$$

Here is an analogue of Theorem 4.

\vskip5pt

{\bf Theorem 5.}
{\it
Let
$C = C^* > 0$. 
The spectrum of the matrix
$A$
belongs to 
${\cal E}_e$ 
if and only if there is a Hermitian positive definite solution
$H = H^* > 0$ 
to~$(8)$.
} 

\vskip5pt

{\bf Proof.}
Assume that all eigenvalues
$\lambda_s$ 
of 
$A$ 
belong to 
${\cal E}_e$. 
Then all eigenvalues
$\mu_r = \bar\lambda_r$ 
of 
$A^*$ 
also belong to 
${\cal E}_e$. 
To prove the unique solvability of~(8) for any matrix
$C$,
it is enough to check that the condition of Krein's theorem
$$
P(\lambda_s,\mu_r)\neq 0
$$
is valid.

Obviously, for~(8), we have
$$
P(\lambda_s,\mu_r) = 1- \left(\frac{1}{2a^2}+\frac{1}{2b^2}\right)
\lambda_s\mu_r - \left(\frac{1}{4a^2}-\frac{1}{4b^2}\right)(\mu^2_r + 
\lambda^2_s)
$$
$$
=1- \left(\frac{1}{2a^2}+\frac{1}{2b^2}\right)
\overline\mu_s\mu_r - \left(\frac{1}{4a^2}-\frac{1}{4b^2}\right)(\mu^2_r + 
\overline\mu^2_s).
$$
Denote
$$
\mu_k = \alpha_k + i\beta_k. 
$$
Taking into account the inequality
$$
\frac{\alpha_k^2}{a^2} + \frac{\beta_k^2}{b^2} > 1, 
$$
it is not difficult to show that the condition of Krein's theorem is fulfilled
and for any matrix
$C$
(8) has a unique solution; moreover,
$$
H = \frac{1}{(2\pi i)^2} \int\limits_{\gamma_A} \int\limits_{\gamma_{A^*}}
\frac{1}{1- \left(\frac{1}{2a^2}+\frac{1}{2b^2}\right)
\lambda\mu - \left(\frac{1}{4a^2}-\frac{1}{4b^2}\right)(\mu^2 + 
\lambda^2)}
$$
$$
\circ (\lambda I - A^*)^{-1} C (\mu I - A)^{-1}
d\lambda\, d\mu,
$$
where the contours
$\gamma_A$,
$\gamma_A{^*}$ 
surround the spectra of the matrices
$A$ 
and
$A^*$,
respectively. Note that 
$H = H^* > 0$
if
$C = C^* > 0$. 

We now prove the theorem in the opposite direction.
Suppose that there is an eigenvalue
$\mu_k$
of 
$A$
such that
$$
\frac{({\rm Re}\,\mu_k)^2}{a^2} + \frac{({\rm Im}\,\mu_k)^2}{b^2} \le 1.
$$
We take an eigenvector
$v_k$
corresponding to
$\mu_k$
and consider the scalar product
$\langle Cv_k, v_k \rangle$. 
By the condition of the theorem,
$\langle Cv_k, v_k \rangle > 0$.

We rewrite this scalar product taking into account that
$H$
is a solution to~(8). We have
$$
- \langle Cv_k, v_k \rangle
= \langle H v_k, v_k \rangle
- \left(\frac{1}{2a^2}+\frac{1}{2b^2}\right) \left\langle HA\,v_k, A\,v_k \right\rangle 
$$
$$ 
- \left(\frac{1}{4a^2}-\frac{1}{4b^2}\right) \left( \left\langle HA^2 \,v_k, v_k \right\rangle 
+ \left\langle H\,v_k, A^2v_k \right\rangle\right).
$$
Taking into account that
$Av_k = \mu_k v_k$,
we get
$$
- \langle Cv_k, v_k \rangle 
= \langle H v_k, v_k \rangle
- \left(\frac{1}{2a^2}+\frac{1}{2b^2}\right) \overline\mu_k
\mu_k \langle H v_k, v_k \rangle
$$
$$
- \left(\frac{1}{4a^2}-\frac{1}{4b^2}\right) (\mu^2_k + \overline \mu^2_k)
\langle Hv_k, v_k \rangle.
$$
Since
$$
\overline\mu_k \mu_k = ({\rm Re}\,\mu_k)^2 + ({\rm Im}\,\mu_k)^2,
$$
$$
\mu^2_k  + \overline\mu^2_k = 2({\rm Re}\,\mu_k)^2-2({\rm Im}\,\mu_k)^2, 
$$
then we obtain
$$
\langle Cv_k, v_k \rangle 
= - \left(1 - \frac{({\rm Re}\,\mu_k)^2}{a^2}-\frac{({\rm Im}\,\mu_k)^2}{b^2}\right) 
\langle Hv_k, v_k \rangle.
$$
Since
$\langle Hv_k, v_k \rangle > 0$,
then by virtue of our assumption
$$
\langle Cv_k, v_k \rangle \le 0. 
$$ 
However, by condition,  
$C > 0$.
Then we have a contradiction.
Therefore, the spectrum of 
$A$
belongs to 
${\cal E}_e$. 

The theorem is proven.

\vskip5pt

We now present two theorems on perturbations of the matrix spectrum.

\vskip5pt

{\bf Theorem 6.}
{\it
Let the spectrum of the matrix
$A$
belong to 
${\cal E}_i$  
and let 
$H = H^* > 0$
be a solution to $(7)$ with
$C = I$. 
If 
$B$
is a matrix such that
$$
\left(\frac{1}{2a^2}+\frac{1}{2b^2}\right) \left(B^*HA + A^*HB + B^*HB\right)
$$
$$
+ \left(\frac{1}{4a^2}-\frac{1}{4b^2}\right)\left(H(AB + BA + B^2) + (AB + BA + B^2)^*H\right) < I,
\eqno(9)
$$
then the spectrum of the perturbed matrix
$A + B$
also belongs to 
${\cal E}_i$.
} 
 
\vskip5pt

{\bf Proof.}
Let
$\lambda_j$
be an eigenvalue of the matrix
$A+B$ 
and let
$v_j$
be the corresponding eigenvector, i.e.
$(A+B)v_j = \lambda_j v_j$.
Consider the scalar product
$$
J_j = \langle [H - \left(\frac{1}{2a^2}+\frac{1}{2b^2}\right) (A+B)^*H(A+B)]v_j,v_j \rangle 
$$
$$ 
- \langle \left(\frac{1}{4a^2}-\frac{1}{4b^2}\right)(H(A+B)^2 + ((A+B)^*)^2H)v_j,v_j \rangle. 
$$
We rewrite it in the following form
$$
J_j = \langle [H - \left(\frac{1}{2a^2} + \frac{1}{2b^2}\right) A^*HA 
- \left(\frac{1}{4a^2}-\frac{1}{4b^2}\right)(HA^2 + (A^*)^2H)]v_j,v_j \rangle
$$
$$
- \langle [\left(\frac{1}{2a^2}+\frac{1}{2b^2}\right) \left(B^*HA + A^*HB + B^*HB\right)]v_j,v_j \rangle
$$
$$
- \langle [\left(\frac{1}{4a^2}-\frac{1}{4b^2}\right)\left(H(AB + BA + B^2) + (AB + BA + B^2)^*H\right)]v_j,v_j \rangle.
$$
Since
$H$
is a solution to~(7) with
$C = I$,
we get
$$
J_j = \langle [I - \left(\frac{1}{2a^2}+\frac{1}{2b^2}\right) \left(B^*HA + A^*HB + B^*HB\right)]v_j,v_j \rangle
$$
$$
- \langle [\left(\frac{1}{4a^2}-\frac{1}{4b^2}\right)\left(H(AB + BA + B^2) + (AB + BA + B^2)^*H\right)]v_j,v_j \rangle.
$$
By~(9) we have
$J_j > 0$. 
On the other hand,
$$
J_j = \langle H v_j, v_j \rangle
- \left(\frac{1}{2a^2}+\frac{1}{2b^2}\right) \left\langle H(A+B)\,v_j, (A+B)\,v_j \right\rangle 
$$
$$ 
- \left(\frac{1}{4a^2}-\frac{1}{4b^2}\right) \left( \left\langle H(A+B)^2 \,v_j, v_j \right\rangle 
+ \left\langle H\,v_j, (A+B)^2v_j \right\rangle\right).
$$
Taking into account the equality
$(A+B)v_j = \lambda_j v_j$,
we obtain
$$
J_j = \langle H v_j, v_j \rangle
- \left(\frac{1}{2a^2}+\frac{1}{2b^2}\right) \overline\lambda_j
\lambda_j \langle H v_j, v_j \rangle
$$
$$
- \left(\frac{1}{4a^2}-\frac{1}{4b^2}\right) (\lambda^2_j + \overline \lambda^2_j)
\langle Hv_j, v_j \rangle
$$
or
$$
J_j = \left(1 - \frac{({\rm Re}\,\lambda_j)^2}{a^2}-\frac{({\rm Im}\,\lambda_j)^2}{b^2}\right) 
\langle Hv_j, v_j \rangle.
$$
Since
$H > 0$, $J_j>0$, 
then
$$
1 - \frac{({\rm Re}\,\lambda_j)^2}{a^2} - \frac{({\rm Im}\,\lambda_j)^2}{b^2} > 0. 
$$
Therefore, due to arbitrariness of 
$\lambda_j$,
we get that the spectrum of the matrix
$A+B$
belongs to 
${\cal E}_i$.

The theorem is proven.

\vskip5pt

{\bf Corollary 1.}
{\it
Let the spectrum of the matrix
$A$ 
belong to 
${\cal E}_i$  
and let
$H = H^* > 0$
be a solution to~$(7)$ with
$C = I$. 
If matrix
$B$
is a matrix such that
$$
\|B\| < \sqrt{\|A\|^2 + \frac{b^2}{\|H\|}} - \|A\|, 
\eqno(10)
$$
then the spectrum of the perturbed matrix
$A + B$
also belongs to 
${\cal E}_i$.
}
 
\vskip5pt

{\bf Proof.}
We consider the scalar product
$$
J(v) = \langle [\left(\frac{1}{2a^2}+\frac{1}{2b^2}\right) \left(B^*HA + A^*HB + B^*HB\right)]v,v \rangle
$$
$$
+ \langle [\left(\frac{1}{4a^2}-\frac{1}{4b^2}\right)\left(H(AB + BA + B^2) + (AB + BA + B^2)^*H\right)]v,v \rangle,
$$
where 
$v \in  C^n$
is an arbitrary vector such that
$\|v\|=1$.
Obviously, for
$a \ge b$
we have
$$
J(v) \le \frac{\|H\|}{b^2} \left(2\|B\|\|A\| + \|B\|^2\right) \|v\|^2. 
$$
Then,
$$
J(v) \le \frac{\|H\|}{b^2} \left(2\|B\|\|A\| + \|B\|^2 - \frac{b^2}{\|H\|}\right) \|v\|^2 + \|v\|^2. 
$$
Since the function
$$
\varphi(t) = t^2 + 2 \|A\| t - \frac{b^2}{\|H\|}
$$
is strictly negative for
$$
-\sqrt{\|A\|^2 + \frac{b^2}{\|H\|}} - \|A\| < t < \sqrt{\|A\|^2 + \frac{b^2}{\|H\|}} - \|A\|, 
$$
then, taking into account~(10), we obtain
$$
J(v) < \|v\|^2 = \langle Iv, v\rangle.
$$
Consequently,
$$
\left(\frac{1}{2a^2}+\frac{1}{2b^2}\right) \left(B^*HA + A^*HB + B^*HB\right)
$$
$$
+ \left(\frac{1}{4a^2}-\frac{1}{4b^2}\right)\left(H(AB + BA + B^2) + (AB + BA + B^2)^*H\right) < I.
$$
By virtue of Theorem~6, the spectrum of the matrix
$A + B$
belongs to 
${\cal E}_i$.

The corollary is proven.

\vskip5pt

The next assertion follows from the above.

\vskip5pt

{\bf Corollary 2.}
{\it
Let the spectrum of the matrix
$A$ 
belong to 
${\cal E}_i$  
and let
$H = H^* > 0$
be a solution to~$(7)$ with
$C = I$. 
If 
$B$
is a matrix such that~$(10)$ is fulfilled, then the matrix equation
$$
H - \left(\frac{1}{2a^2} + \frac{1}{2b^2}\right) (A + B)^*H(A + B)  
$$
$$
- \left(\frac{1}{4a^2} - \frac{1}{4b^2}\right)(H(A + B)^2 + (A^* + B^*)^2H) = C 
$$
is uniquely solvable for any right-hand side
$C$.
}
 
\vskip5pt

{\bf Theorem 7.}
{\it 
Let the spectrum of the matrix
$A$ 
belong to 
${\cal E}_e$  
and let 
$H = H^* > 0$
be a solution to~$(8)$ with
$C = I$. 
If 
$B$
is a matrix such that
$$
-I < \left(\frac{1}{2a^2}+\frac{1}{2b^2}\right) \left(B^*HA + A^*HB + B^*HB\right)
$$
$$
+ \left(\frac{1}{4a^2}-\frac{1}{4b^2}\right)\left(H(AB + BA + B^2) + (AB + BA + B^2)^*H\right),  
\eqno(11)
$$
then the spectrum of the perturbed matrix
$A + B$
also belongs to 
${\cal E}_e$.
}
 
\vskip5pt

{\bf Proof.}
We conduct reasoning according to the scheme of the proof of Theorem 6.
Let
$\lambda_j$
be an eigenvalue of the matrix
$A+B$ 
and let
$v_j$
be the corresponding eigenvector. Consider the scalar product
$$
J_j = \langle [H - \left(\frac{1}{2a^2}+\frac{1}{2b^2}\right) (A+B)^*H(A+B)]v_j,v_j \rangle 
$$
$$ 
- \langle \left(\frac{1}{4a^2}-\frac{1}{4b^2}\right)(H(A+B)^2 + ((A+B)^*)^2H)v_j,v_j \rangle. 
$$
Since
$H$
is a solution to~(8) with
$C = I$,
then by analogy with the previous one, 
$J_j$ 
can be rewritten as follows
$$
J_j = \langle [-I - \left(\frac{1}{2a^2}+\frac{1}{2b^2}\right) \left(B^*HA + A^*HB + B^*HB\right)]v_j,v_j \rangle
$$
$$
- \langle [\left(\frac{1}{4a^2}-\frac{1}{4b^2}\right)\left(H(AB + BA + B^2) + (AB + BA + B^2)^*H\right)]v_j,v_j \rangle.
$$
By~(11) we have
$J_j < 0$.
On the other hand,
$$
J_j = \langle H v_j, v_j \rangle
- \left(\frac{1}{2a^2}+\frac{1}{2b^2}\right) \left\langle H(A+B)\,v_j, (A+B)\,v_j \right\rangle 
$$
$$ 
- \left(\frac{1}{4a^2}-\frac{1}{4b^2}\right) \left( \left\langle H(A+B)^2 \,v_j, v_j \right\rangle 
+ \left\langle H\,v_j, (A+B)^2v_j \right\rangle\right) 
$$
or
$$
J_j = \left(1 - \frac{({\rm Re}\,\lambda_j)^2}{a^2}-\frac{({\rm Im}\,\lambda_j)^2}{b^2}\right) 
\langle Hv_j, v_j \rangle.
$$
Since
$H > 0$, $J_j < 0$, 
then
$$
1 - \frac{({\rm Re}\,\lambda_j)^2}{a^2} - \frac{({\rm Im}\,\lambda_j)^2}{b^2} < 0. 
$$ 
Therefore, due to arbitrariness of 
$\lambda_j$,
we get that the spectrum of the matrix
$A+B$
belongs to 
${\cal E}_e$.

The theorem is proven.

\vskip5pt

{\bf Corollary 1.}
{\it
Let the spectrum of the matrix
$A$ 
belong to 
${\cal E}_e$  
and let
$H = H^* > 0$
be a solution to~$(8)$ with
$C = I$. 
If 
$B$
is a matrix such that
$$
\|B\| 
< 2a^2(a^2 - b^2)^{-1} \left(\sqrt{\|A\|^2 + \left(\frac{a^2 - b^2}{2a^2}\right)\frac{b^2}{\|H\|}} 
- \|A\|\right), 
\eqno(12)
$$
then the spectrum of the perturbed matrix
$A + B$
also belongs to 
${\cal E}_e$.
}
 
\vskip5pt

{\bf Proof.}
We consider the scalar product
$$
J(v) = -\langle [\left(\frac{1}{2a^2}+\frac{1}{2b^2}\right) \left(B^*HA + A^*HB + B^*HB\right)]v,v \rangle
$$
$$
- \langle [\left(\frac{1}{4a^2}-\frac{1}{4b^2}\right)\left(H(AB + BA + B^2) + (AB + BA + B^2)^*H\right)]v,v \rangle,
$$
where 
$v \in  C^n$
is an arbitrary vector such that 
$\|v\|=1$,
Obviously, for
$a \ge b$
we have
$$
J(v) \le \frac{\|H\|}{b^2} \left(2\|B\|\|A\| + \left(\frac{a^2 - b^2}{2a^2}\right) \|B\|^2\right) \|v\|^2. 
$$
Then,
$$
J(v) \le \frac{\|H\|}{b^2} \left(2\|B\|\|A\| + \left(\frac{a^2 - b^2}{2a^2}\right) \|B\|^2 
- \frac{b^2}{\|H\|}\right) \|v\|^2 + \|v\|^2. 
$$
Since the function
$$
\varphi(t) = \left(\frac{a^2 - b^2}{2a^2}\right) t^2 + 2 \|A\| t - \frac{b^2}{\|H\|}
$$ 
for
$$
0 \le t < 2a^2(a^2 - b^2)^{-1} \left(\sqrt{\|A\|^2 + \left(\frac{a^2 - b^2}{2a^2}\right)\frac{b^2}{\|H\|}} 
- \|A\|\right) 
$$
is strictly negative, then, taking into account~(12), we obtain 
$$
J(v) < \|v\|^2 = \langle Iv, v\rangle.
$$
Consequently,
$$
- \left(\frac{1}{2a^2}+\frac{1}{2b^2}\right) \left(B^*HA + A^*HB + B^*HB\right)
$$
$$
- \left(\frac{1}{4a^2}-\frac{1}{4b^2}\right)\left(H(AB + BA + B^2) + (AB + BA + B^2)^*H\right) < I.
$$
By virtue of Theorem~7, the spectrum of the matrix
$A + B$
belongs to 
${\cal E}_e$.

The corollary is proven.

\vskip5pt

The next assertion follows from the above.

\vskip5pt

{\bf Corollary 2.}
{\it
Let the spectrum of the matrix
$A$ 
belong to 
${\cal E}_e$  
and let
$H = H^* > 0$
be a solution to~$(8)$ with
$C = I$. 
If 
$B$
is a matrix such that~$(12)$ is fulfilled, then the matrix equation
$$
H - \left(\frac{1}{2a^2} + \frac{1}{2b^2}\right) (A + B)^*H(A + B)  
$$
$$
- \left(\frac{1}{4a^2} - \frac{1}{4b^2}\right)(H(A + B)^2 + (A^* + B^*)^2H) = C 
$$
is uniquely solvable for any right-hand side
$C$.
}

\section{Localization of the matrix spectrum with respect to a parabola}

The following problem was also solved in~\cite{M},~\cite{D09}.

\vskip5pt

{\bf Problem 2.} Belonging of the spectrum of the matrix
$A$
to the domain bounded by a parabola
$$
{\cal P}_i =
\left\{\lambda : \ ({\rm Im}\,\lambda)^2 < 2p {\rm Re}\, \lambda \right\}, 
\quad p > 0,
$$
is equivalent to the solvability of the matrix equation
$$
HA + A^*H - \frac{1}{2p}A^*HA + \frac{1}{4p}(HA^2+(A^*)^2H) = I
\eqno(13)
$$
in the class of Hermitian positive definite matrices
$H = H^*>0$. 

\vskip5pt

In particular, the following theorem was proved in~\cite{M},~\cite{D09}.

\vskip5pt

{\bf Theorem 8.}
{\it 
The spectrum of the matrix
$A$
belongs to 
${\cal P}_i$ 
if and only if the matrix equation
$$
HA + A^*H - \frac{1}{2p}A^*HA + \frac{1}{4p}(HA^2+(A^*)^2H) = C
$$
with 
$C = C^* > 0$ 
has a Hermitian positive definite solution
$H = H^* > 0$.  
}

\vskip5pt

Note that (13) has the form (6) with 
$$
N=2, \quad n=m, \quad B=A^*, \quad Y=I,
$$
$$
a_{00}=a_{12}=a_{21}=a_{22}=0, \quad a_{01}=a_{10}=1, 
$$
$$
a_{11}=-\frac{1}{2p}, \quad 
a_{02}=a_{20}=\frac{1}{4p}.
$$

Consider the problem on the belonging of the matrix spectrum to a set
lying outside the domain bounded by a parabola
$$
\left\{\lambda : \ ({\rm Im}\,\lambda)^2 = 2p {\rm Re}\,\lambda\right\}, \quad p > 0, 
$$
i.e. the problem on the belonging to the set
$$
{\cal P}_e =
\left\{\lambda : \ ({\rm Im}\,\lambda)^2 > 2p {\rm Re}\, \lambda \right\}. 
$$

By analogy with the problem of the belonging of the matrix spectrum to 
${\cal P}_i$   
consider the matrix equation of the form~(13) with a modified right-hand side
$$
HA + A^*H - \frac{1}{2p}A^*HA + \frac{1}{4p}(HA^2+(A^*)^2H) = - C.
\eqno(14)
$$

The following theorem was proved in~\cite{D09},~\cite{DP}.

\vskip5pt

{\bf Theorem 9.}
{\it 
Suppose that there is a Hermitian positive definite solution
$H = H^* > 0$ 
to~$(14)$
with 
$C = C^* > 0$. Then the spectrum of the matrix
$A$
belongs to 
${\cal P}_e$. 
}

\vskip5pt

Note that, if the spectrum of 
$A$
lies outside the closure of 
${\cal P}_i$, 
then one can check feasibility of the condition 
$$
P(\lambda_s,\mu_r) \neq 0
$$
from Theorem~3 and show the solvability of~(14) for any matrix
$C$. 
Moreover, if
$C = C^* > 0$, 
then
$H = H^* > 0$, 
i.e., it is possible to formulate a criterion for the belonging
of the matrix spectrum to the domain lying outside the closure of
${\cal P}_i$.

By analogy with the previous section,
we present two theorems on perturbations of the matrix spectrum.

\vskip5pt

{\bf Theorem 10.}
{\it 
Let the spectrum of the matrix
$A$
belong to 
${\cal P}_i$  
and let 
$H = H^* > 0$
be a solution to~$(13)$. If 
$B$
is a matrix such that
$$
HB + B^*H - \frac{1}{2p}(A^*HB + B^*HA + B^*HB)  
$$
$$
+ \frac{1}{4p}(H(AB + BA + B^2) + (AB + BA + B^2)^*H) > -I,
\eqno(15)
$$
then the spectrum of the perturbed matrix
$A + B$
also belongs to 
${\cal P}_i$.
}
 
\vskip5pt

{\bf Proof.}
Let
$\lambda_j$
be an eigenvalue of the matrix
$A+B$ 
and let 
$v_j$
be the corresponding eigenvector, i.e.
$$
(A + B)v_j = \lambda_j v_j. 
$$
Consider the scalar product
$$
J_j = \langle H(A+B)v_j,v_j \rangle + \langle (A+B)^*Hv_j,v_j \rangle 
- \langle \frac{1}{2p}(A+B)^*H(A+B)v_j,v_j \rangle 
$$
$$ 
+ \langle \left(\frac{1}{4p}(H(A+B)^2 + ((A+B)^*)^2H)\right)v_j,v_j \rangle. 
$$
Since
$H$
is a solution to~(13), we get
$$
J_j = \langle (I + HB + B^*H)v_j,v_j \rangle 
- \langle \left(\frac{1}{2p}(A^*HB + B^*HA + B^*HB)\right)v_j,v_j \rangle 
$$
$$
+ \langle \left(\frac{1}{4p} (H(AB + BA + B^2) + (AB + BA + B^2)^*H)\right)v_j,v_j \rangle.
$$
By~(15) we have
$J_j > 0$.
On the other hand,
$$
J_j = \langle H(A+B)v_j,v_j \rangle + \langle Hv_j, (A+B)v_j \rangle 
- \langle \frac{1}{2p}H(A+B)v_j, (A+B)v_j \rangle 
$$
$$ 
+ \frac{1}{4p}(\langle H(A+B)^2 v_j, v_j \rangle + \langle Hv_j, (A+B)^2 v_j \rangle). 
$$
Taking into account the equality
$(A+B)v_j = \lambda_j v_j$,
we obtain
$$
J_j = (\lambda_j + \bar{\lambda_j})\langle H v_j, v_j \rangle
- \frac{1}{2p}\lambda_j \overline\lambda_j \langle H v_j, v_j \rangle  
+ \frac{1}{4p}(\lambda^2_j + \overline \lambda^2_j) \langle Hv_j, v_j \rangle
$$
or
$$
J_j = \left(2{\rm Re}\,\lambda_j - \frac{1}{p} ({\rm Im}\,\lambda_j)^2\right) 
\langle Hv_j, v_j \rangle.
$$
Since
$H > 0$, $J_j>0$,
then
$$
2{\rm Re}\,\lambda_j - \frac{1}{p} ({\rm Im}\,\lambda_j)^2 > 0. 
$$ 
Therefore, due to arbitrariness of 
$\lambda_j$,
we get that the spectrum of the matrix
$A+B$
belongs to 
${\cal P}_i$.

The theorem is proven.

\vskip5pt

{\bf Corollary.}
{\it
Let the spectrum of the matrix
$A$ 
belong to 
${\cal P}_i$  
and let
$H = H^* > 0$
be a solution to~$(13)$. If 
$B$
is a matrix such that~$(15)$ is fulfilled, then the matrix equation
$$
H(A + B) + (A + B)^*H - \frac{1}{2p}(A + B)^*H(A + B) 
$$
$$
+ \frac{1}{4p}(H(A + B)^2 + (A^* + B^*)^2H) = C 
\eqno(16)
$$
is uniquely solvable for any right-hand side
$C$.
}

\vskip5pt

{\bf Theorem 11.}
{\it 
Let the spectrum of the matrix
$A$ 
belong to 
${\cal P}_e$  
and let
$H = H^* > 0$
be a solution to~$(14)$ with
$C = I$. 
If 
$B$
is a matrix such that
$$
HB + B^*H - \frac{1}{2p}(A^*HB + B^*HA + B^*HB)  
$$
$$
+ \frac{1}{4p}(H(AB + BA + B^2) + (AB + BA + B^2)^*H) < I,
\eqno(17)
$$
then the spectrum of the perturbed matrix
$A + B$
also belongs to 
${\cal P}_e$.
}
 
\vskip5pt

{\bf Proof.}
Let
$\lambda_j$
be an eigenvalue of the matrix
$A+B$ 
and let
$v_j$
be the corresponding eigenvector. Consider the scalar product
$$
J_j = \langle H(A+B)v_j,v_j \rangle + \langle (A+B)^*Hv_j,v_j \rangle 
- \langle \frac{1}{2p}(A+B)^*H(A+B)v_j,v_j \rangle 
$$
$$ 
+ \langle \left(\frac{1}{4p}(H(A+B)^2 + ((A+B)^*)^2H)\right)v_j,v_j \rangle. 
$$
Since
$H$
is a solution to~(14) with
$C = I$, 
we get
$$
J_j = \langle (-I + HB + B^*H)v_j,v_j \rangle 
- \langle \left(\frac{1}{2p}(A^*HB + B^*HA + B^*HB)\right)v_j,v_j \rangle 
$$
$$
+ \langle \left(\frac{1}{4p} (H(AB + BA + B^2) + (AB + BA + B^2)^*H)\right)v_j,v_j \rangle.
$$
By~(17) we have
$J_j < 0$. 
On the other hand,
$$
J_j = \langle H(A+B)v_j,v_j \rangle + \langle Hv_j, (A+B)v_j \rangle 
- \frac{1}{2p} \langle H(A+B)v_j, (A+B)v_j \rangle 
$$
$$ 
+ \frac{1}{4p}(\langle H(A+B)^2 v_j, v_j \rangle + \langle Hv_j, (A+B)^2 v_j \rangle). 
$$
Taking into account the equality
$(A + B)v_j = \lambda_j v_j$,
we obtain
$$
J_j = (\lambda_j + \bar{\lambda_j})\langle H v_j, v_j \rangle
- \frac{1}{2p}\lambda_j \overline\lambda_j \langle H v_j, v_j \rangle  
+ \frac{1}{4p}(\lambda^2_j + \overline \lambda^2_j) \langle Hv_j, v_j \rangle
$$
or
$$
J_j = \left(2{\rm Re}\,\lambda_j - \frac{1}{p} ({\rm Im}\,\lambda_j)^2\right) 
\langle Hv_j, v_j \rangle.
$$
Since
$H > 0$, $J_j < 0$, 
then
$$
2{\rm Re}\,\lambda_j - \frac{1}{p} ({\rm Im}\,\lambda_j)^2 < 0. 
$$ 
Therefore, due to arbitrariness of 
$\lambda_j$
we get that the spectrum of the matrix
$A + B$ 
belongs to 
${\cal P}_e$.

The theorem is proven.

\vskip5pt

{\bf Corollary.}
{\it
Let the spectrum of the matrix
$A$ 
belong to 
${\cal P}_e$  
and let
$H = H^* > 0$
be a solution to~$(14)$ with
$C = I$. 
If 
$B$
is a matrix such that~$(17)$ is fulfilled,
then $(16)$ is uniquely solvable for any right-hand side
$C$.
}

\vskip5pt

Some results on the localization of the matrix spectrum
can also be found in~\cite{SaTo}, \cite{BiBlPo}, \cite{Bi}, \cite{TrDo}.

\end{document}